\documentclass[letterpaper, 12pt]{amsart}

\usepackage[T1]{fontenc}

\usepackage[margin=1in]{geometry}
 
\usepackage{tikz}

\usetikzlibrary{shadings}

\usepackage{amscd}

\usepackage{amsmath}

\usepackage{amsfonts}

\usepackage{amssymb}

\usepackage{amsthm}

\usepackage{bbold}

\usepackage{mathtools}

\usepackage{arydshln}

\usepackage{fancyhdr}

\usepackage{verbatim}

\usepackage[all]{xy}

\usepackage{enumitem}

\usepackage{color}

\usepackage[colorlinks=true, linkcolor=teal, citecolor=-teal,
pagebackref=true]{hyperref}

\usepackage[reftex]{theoremref}

\usepackage{nameref}

\usepackage{charter}
\usepackage{newtxmath}

\theoremstyle{plain}
\newtheorem*{theorem*}{Theorem}
\newtheorem{theorem}{Theorem}
\newtheorem{lemma}[theorem]{Lemma}
\newtheorem{corollary}[theorem]{Corollary}
\newtheorem{proposition}[theorem]{Proposition}

\newtheorem*{claim*}{Claim}

\theoremstyle{definition}

\theoremstyle{remark}
\newtheorem*{remark*}{Remark}

\newcommand{\ba}{\mathbf a}

\newcommand{\RR}{\mathbb{R}}

\newcommand{\QQ}{\mathbb{Q}}

\newcommand{\NN}{\mathbb{N}}

\newcommand{\ZZ}{\mathbb{Z}}

\newcommand{\calH}{\mathcal{H}}

\newcommand{\bone}{\mathbf{1}}

\renewcommand{\leq}{\leqslant} \renewcommand{\geq}{\geqslant}

\DeclarePairedDelimiter{\abs}{\lvert}{\rvert}
\DeclarePairedDelimiter{\norm}{\lVert}{\rVert}

\DeclarePairedDelimiter{\set}{\lbrace}{\rbrace}
\DeclarePairedDelimiter{\parens}{\lparen}{\rparen}

\DeclareMathOperator{\meas}{\lambda}

\def\eps{{\varepsilon}}

\def\1int{{[0,1]}}

\title{Twisted approximation with restricted denominators}

\author{Manuel Hauke \and Felipe A. Ram\'irez}

\date{}

\address{Institute of Analysis and Number Theory, University of Technology Graz,
Graz, Austria}

\email{hauke@math.tugraz.at}

\address{Department of Mathematics and Computer Science, Wesleyan University, Connecticut, USA}

\email{framirez@wesleyan.edu}



\begin{document}

\begin{abstract}
Given an increasing integer sequence $(a_n)$, a real number $\alpha$, and a sequence $\psi(n)$, we study the set $W$ of real numbers $\gamma$ for which $a_n\alpha - \gamma$ is a distance less than $\psi(n)$ away from an integer. This is often referred to as twisted Diophantine approximation, in this case with denominators restricted to the given sequence $(a_n)$. Our main results are about the size of $W$, and they hold for almost every $\alpha$, with respect to a measure of positive Fourier dimension, for example Lebesgue measure. Our results extend recent work of Kristensen and Persson, and answer questions that they posed.
\end{abstract}

\maketitle

\section{Introduction and results}

Given an increasing integer sequence $\ba=(a_n)$ and function
$\psi:\NN\to\RR$, we define
\begin{equation*}
    W(\psi, \ba) = \set*{(\alpha, \gamma) \in[0,1]^2 : \norm{a_n\alpha - \gamma} < \psi(n) \quad \textrm{i.o.}}
\end{equation*}
where $\norm{\cdot}$ denotes distance to the nearest integer, and ``ì.o.'' stands for ``infinitely often''. Sets of the form $W(\psi, \ba)$ are central to the study of metric Diophantine approximation. When $\ba = (n)_{n\in\NN}$, the horizontal fibers of $W(\psi, \ba)$ through a fixed $\gamma$ are the focus of \emph{inhomogeneous approximation}. Study of the vertical fibers through a fixed $\alpha$ is often referred to as \emph{twisted approximation}. For general sequences $\ba=(a_n)$, one speaks of approximations whose denominators have been restricted to $\ba$. 

In this paper we are primarily concerned with vertical fibers of $W(\psi,\ba)$. For each $\alpha$, define
\begin{equation*}
    W(\psi,\ba,\alpha) = \set*{\gamma\in [0,1] : (\alpha, \gamma) \in W(\psi, \ba)}.
\end{equation*}
These have been of great interest since the 1950s. In~\cite{Kurzweil}, Kurzweil used these sets to give an alternate characterization of the set of badly approximable numbers, inspiring many follow-ups decades later such as ~\cite{ChaikaKurzweil, FayadMSTP,
  FK16, Harraptwisted, HK25, Kim, liwang, NotTooBad, bad_characterization, SimmonsKurzweil,
  TsengSTPs}. Part of the importance of the sets $W(\psi,\ba,\alpha)$ is due to their connection to dynamics---they are the sets of $\gamma$ whose orbits under a rotation $\alpha$ enter a shrinking target around $0\in \RR/\ZZ$ infinitely many times. Their study can be seen as a refinement of the classical Poincar{\' e} recurrence theorem from ergodic theory, which says that under an ergodic transformation, almost every point in the phase space will enter a fixed (non-shrinking) target infinitely often. This dynamical view is especially apparent in~\cite{ChaikaKurzweil, FayadMSTP, TsengSTPs}. In this view, restriction to $\ba$ amounts to requiring that the shrinking target is hit at the prescribed times $(a_n)\subset\NN$.

  The study of restricting to the behaviour of subsequences arises from various motivations:
  For instance, a naive approach to Littlewood's Conjecture (and various inhomogeneous variants) is to consider a sequence $(a_n)$ of good approximations (e.g. convergent denominators) to a fixed irrational $\beta$, in order to obtain metric or dimensional results on the sets of $\alpha$ such that $(a_n\alpha)$ visits shrinking targets around $0$, or some inhomogeneous target $\gamma$ (in the case of the inhomogeneous variants). Since these results are most interesting when $\alpha$ is not chosen with respect to Lebesgue measure, but to a measure $\mu$ that is supported on the ``thinner'' set of badly approximable $\alpha$, various results have been established by using different measures $\mu$. For results in this direction, we refer the reader to e.g. \cite{CT24,CZ21,musings,PV,S25}. 

  A further motivation for restricting attention to the behavior of rotation orbits along subsequences arises from the Berry--Tabor Conjecture in quantum chaos theory, having to do with statistics for the spacings between energy levels of quantum systems. The conjecture inspired the seminal works of Rudnick,Sarnak and Zaharescu \cite{rud_sar,rsz}, who examined the dynamics of $(n^2\alpha)$ for almost all $\alpha$, conjecturing certain fine-scale statistics to hold for all badly approximable $\alpha$ (resp. all Diophantine $\alpha$), and showing some of them to hold for Lebesgue-almost all $\alpha$. Their result has been generalized to various different sequences of the form $(a_n\alpha)$. Notably, in~\cite{all,rough} statements on $\alpha$ are made that go beyond the Lebesgue measure theory.

  \

  Recently, Kristensen and Persson~\cite{KP25} have undertaken the
  study of the sets $W(\psi,\ba,\alpha)$ with the aim of finding
  measure-theoretic and Hausdorff dimension results that hold for
  $\mu$-almost every $\alpha$, where $\mu$ is either the Lebesgue
  measure or some other measure of positive Fourier dimension. Among
  other results, they prove:
\begin{itemize}
\item If $\mu$ is a measure on $[0,1]$ with positive Fourier
  dimension, and $\ba$ satisfies a separation property of the form
  $\abs{a_m-a_n} \gg \abs{n-m}^{1/\tau + \eps}$, then for all
  $\sigma < 1$ and $\mu$-almost every $\alpha$, the set
  $W(\psi_\sigma, \ba, \alpha)$ has positive Lebesgue
  measure. (See~\cite[Theorem~3.1]{KP25}.)

\item If $\abs{a_n - a_m} \gg \abs{n-m}^{1 + \eps}$, then for all
  $\sigma < 1$ and Lebesgue-almost all $\alpha\in[0,1]$, the set
  $W(\psi_\sigma, \ba, \alpha)$ has full Lebesgue measure, where
  $\psi_\sigma(n) = n^{-\sigma}$. (See~\cite[Theorem~3.2]{KP25}.)

\item Let $(a_n)$ be an increasing sequence of integers. For any
  $\sigma \geq 3$ and Lebesgue almost every $\alpha$, we have
  $\dim_H(W(\psi_{\sigma},\ba,\alpha)) = \frac{1}{\sigma}$
  (see~\cite[Corollary~3.7]{KP25}).
\end{itemize}
This paper is inspired by those statements. In our main results, we
strengthen the aforementioned Kristensen--Persson theorems in various
aspects. Specifically:
\begin{itemize}
    \item Kristensen and Persson expressed a strong suspicion that the conclusion of~\cite[Theorem~3.1]{KP25} can be strengthened from positive to full measure. In \th\ref{muaefiber}, we confirm this suspicion.

    \item In \th\ref{lebaefiber}, we remove the separation condition
      from~\cite[Theorem~3.2]{KP25} and extend the result so that it
      holds for general functions $\psi(n)$.

    \item Kristensen and Persson remark that it would be natural to
      suspect that~\cite[Corollary~3.7]{KP25} can be extended to all
      $\sigma >1$. In \th\ref{jarnik}, we prove a Jarnik-style theorem
      for the sets $W(\psi,\ba, \alpha)$ which has as corollary the
      expected extension. (See \th\ref{jarnikbesicovitch}.)
\end{itemize}
More precisely, our theorems are as follows.

\begin{theorem}\th\label{muaefiber}
  If $\mu$ is a probability measure on $[0,1]$ with positive Fourier
  dimension with decay $\abs{\hat\mu(\xi)} = O(\abs{\xi}^{-\tau})$, $\ba=(a_n)$ is an integer sequence with
  $\abs{a_m - a_n}>c\abs{m-n}^{\frac{1}{\tau} + \eps}$ for some
  $c, \eps>0$ and $\psi(n) = n^{-\sigma}$ for some $\sigma < 1$, then
  $\meas(W(\psi,\ba,\alpha)) = 1$ for $\mu$-almost every $\alpha$.
\end{theorem}

\begin{theorem}\th\label{lebaefiber}
  If $\ba = (a_n)$ is an increasing integer sequence, and $\psi:\NN\to\RR_{\geq 0}$ 
  such that $\sum\psi(n)=\infty$, then $\meas(W(\psi,\ba,\alpha)) = 1$
  for $\meas$-almost every $\alpha$, where $\meas$ denotes the Lebesgue
  measure.
\end{theorem}

\begin{theorem}\th\label{jarnik}
  Let $\ba = (a_n)$ be an increasing integer sequence and
  $\psi:\NN\to\RR_{\geq 0}$. Let $f$ be a dimension function such that
  $x^{-1} f(x)$ is monotonic. Then 
\begin{equation*}
    \calH^f(W(\psi, \ba, \alpha)) = 
    \begin{cases}
        0 &\textrm{if } \sum f(\psi(n)) < \infty, \\
        \calH^f([0,1]) &\textrm{if } \sum f(\psi(n)) = \infty
    \end{cases}
\end{equation*}
for $\meas$-almost every $\alpha$.
\end{theorem}

\begin{remark*} 
In Theorems \ref{lebaefiber} and \ref{jarnik} we do not need any regularity assumption on $\psi$, i.e. we can even allow arbitrary, non-monotonic functions $\psi: \mathbb{N} \to \mathbb{R}_{\geq 0}$.
\end{remark*}

\begin{corollary}\th\label{jarnikbesicovitch}
  Let $\ba = (a_n)$ be an increasing integer sequence. For $\meas$-almost every $\alpha\in[0,1]$, the following holds: For all $\sigma\geq 1$ we have
\begin{equation*}
    \dim_H(W(\psi_\sigma, \ba, \alpha)) = \frac{1}{\sigma}
\end{equation*}
where $\psi_\sigma (n) = n^{-\sigma}$.
\end{corollary}

\subsection*{Acknowledgements}
A part of this project was carried out while MH and FAR were on research stay at the University of York. The authors would like to thank the whole number theory group of the University of York for their hospitality. This research was
funded in whole or in part by the Austrian Science Fund (FWF)
10.55776/ESP5134624.

\section{Notation and standard definitions}
For a Radon measure $\mu$ on $\mathbb{R}$, we define the Fourier dimension as 
\[\dim_F(\mu) := \sup\left\{s \in [0,1]: \abs{\hat\mu(\xi)} = O(\abs{\xi}^{-s/2})\right\},\]
where the Fourier transform $\hat{\mu}$ is defined by
$\hat{\mu}(\xi) = \int e^{2 \pi i x \xi} d\mu(x)$.
We write $\lambda$ for the $1$-dimensional Lebesgue measure and $B(x_0,r)$
 for the ball with center $x_0$ and radius $r$. We use the Vinagradov notations $\ll,\gg$, with $f\ll g$ denoting $f \leq C g$ for some absolute constant $C$. We write $\lVert \cdot \rVert$ for the distance to the nearest integer.

For a continuous nondecreasing function $f:\RR^+ \to \RR^+$ such that $f(r)\to 0$ as $r\to 0$, we denote the corresponding $f$-Hausdorff measure by $\calH^f$, and call $f$ a \emph{dimension function}. Given a ball $B:=B(x,r)$, we denote $B^f:=B(x, f(r))$. When $f(x) = x^s$ for some $s > 0$, we use the notation $B^s = B^f$. In particular, $B^1 = B$. 

\section{Prerequisites}

We make use of Beresnevich and Velani's mass transference principle~\cite[Theorem 2]{MTP}:

\begin{theorem}[Mass Transference Principle]\th\label{MTPthm}
  Let $\{B_i\}_{i\in \mathbb{N}}$ be a sequence of balls in
  $\mathbb{R}^k$ with $r(B_i) \to 0$ as $i \to \infty$. Let $f$ be a
  dimension function such that $x^{-k} f(x)$ is monotonic and suppose
  that for any ball $B$ in $\mathbb{R}^k$
\[\mathcal{H}^k(B \cap \limsup B_i^{f}) = \mathcal{H}^k(B).\]
Then, for any ball $B$ in $\mathbb{R}^k$
\[\calH^f(B \cap \limsup B_i^k) = \calH^f(B).\]
\end{theorem}

Furthermore, we employ a tool that is handy to move from positive measure to full measure:

\begin{proposition}[Beresnevich--Dickinson--Velani,~{\cite[Lemma 6]{BDV}}]\th\label{BDVdensitylemma}
  Let $(X,d)$ be a metric space with a finite measure $\mu$ such that
  every open set is $\mu$-measurable. Let $A$ be a Borel subset of $X$
  and let $f:\RR_+\to\RR_+$ be an increasing function with $f(x)\to 0$ as
  $x\to 0$. If for every open set $U\subset X$ we have
  \begin{equation*}
    \mu(A\cap U) \geq f(\mu(U)),
  \end{equation*}
  then $\mu(A) = \mu(X)$.
\end{proposition}

Finally, we employ a result of Haynes--Jensen--Kristensen that generalizes Weyl's famous equidistribution result from the Lebesgue measure to arbitrary measures with sufficient Fourier decay:
\begin{proposition}[Haynes--Jensen--Kristensen,~{\cite[Corollary 7]{musings}}]\label{prop_equi_mu}
Let $\ba =(a_n)_{n \in \mathbb{N}}$ be an increasing integer sequence and let $\mu$ be a probability measure on $[0,1]$ with Fourier decay $\abs{\hat\mu(\xi)} = O(\abs{\xi}^{-\tau})$. Then for $\mu$-almost every $\alpha$, the sequence
$(a_n\alpha)$ is equidistributed mod $1$.
\end{proposition}

\section{Proof of Theorem~\ref{muaefiber}}

Following Kristensen and Persson \cite[Pages~24--25]{KP25}, we
consider a fastly increasing sequence $(n_j)_{j \in \mathbb{N}}$. For
technical reasons, we may (and will) assume that for all
$j \in \mathbb{N}$, $n_j = 2^{k_j}$ for some $k_j \in \mathbb{N}$.

We define
\begin{align*}
    B_k(\alpha) &= B(a_k\alpha,k^{-\sigma}), \quad &S_j = \sum_{n_j \leq k < n_{j+1}}\lambda(B_k(\alpha)),\\
    C_{m,n}(\alpha) &= \sum_{m \leq k,l < n}\lambda(B_k(\alpha) \cap B_l(\alpha)),\quad
    &C_j = \int C_{n_j,n_{j+1}}(\alpha)\, d\mu(\alpha),
\end{align*}
noting that $S_j$ is independent of $\alpha$.  Since $(n_j)$ is
allowed to grow arbitrarily fast, we may assume that $S_j \geq 1$ for
all $j$.  Following Kristensen and Persson \cite[Eqs.~(4.3)
and~(4.4)]{KP25} verbatim, we obtain
\[C_j \leq C S_j^2\]
for some absolute constant $C > 0$.
An application of the Markov inequality implies that for every $p > 1$, the set
\[G(p) := \{\alpha: C_{n_j,n_{j+1}}(\alpha) \leq C\cdot p \cdot S_j^2 \text{ for i.m. }j \}\]
satisfies
$\mu(G(p)) \geq 1 - 1/p$. 

By Proposition \ref{prop_equi_mu}, we obtain the existence of a set 
$S \subseteq [0,1]$ with $\mu(S) = 1$ such that
$(a_n\alpha)$ is uniformly distributed for every $\alpha \in S$.

Thus we see that $G'(p) := G(p) \cap S$ also satisfies
$\mu(G'(p)) \geq 1 - 1/p$. We now claim that for every $\alpha \in G'(p)$, we have 
\begin{equation}\label{limsup_ak}\lambda(\limsup_{k \to \infty} B_k(\alpha)) = 1.\end{equation}
For fixed $\alpha \in G'(p)$, let $(j_i)_{i \in \mathbb{N}}$ be the sequence along which we have 
$C_{n_{j_i},n_{{j_i}+1}}(\alpha) \leq C\cdot p \cdot S_{j_i}^2$.
Since $S_j \geq 1$, we observe that $\sum_{i} S_{j_i} = \infty$. We will actually show
\[\lambda\bigg(\limsup_{\substack{k \to \infty\\ k \in \cup_{i}[n_{j_i},n_{j_i+1})}} B_k(\alpha)\bigg) = 1,\]
which obviously suffices to prove \eqref{limsup_ak}.
We now fix one such $j = j_i$, and fix an arbitrary open set $U$. We want to prove that there is an absolute constant $C$ such that for any $j$ large enough, 

\begin{equation}\label{equid}\sum_{n_j \leq k < n_{j+1}}\lambda(B_k) \leq \frac{C}{\lambda(U)}\sum_{n_j \leq k < n_{j+1}}\lambda(B_k\cap U).\end{equation}

We may assume without loss of generality that $U$ is an open interval $B(u,r)$: Since $U \subseteq [0,1)$ is a union of countably many open intervals, a collection of finitely many such intervals provides a subset of $U$ with the sum of its measures being at least $\lambda(U)/2$. By changing the constant $C$, this is enough to prove \eqref{equid} for $U = B(u,r)$. Now we replace $B_k(\alpha) = B(a_k\alpha,k^{-\sigma})$
by $B_k'(\alpha) = B(a_k\alpha,({2^{n+1}})^{-\sigma})$
where $2^n \leq k < 2^{n+1}$, and observe that
\[B_k'(\alpha) \subseteq B_k(\alpha), \quad \lambda(B_k'(\alpha)) \gg \lambda(B_k(\alpha)).
\]
Writing $U' = B(u,r/2)$, we get for $n$ large enough that
\[\sum_{2^n \leq k \leq 2 ^{n+1}}
\lambda(B_k(\alpha) \cap U)
\geq ({2^{n+1}})^{-\sigma}\#\{2^n \leq k \leq 2^{n+1}: a_k\alpha \in U'\}.
\]
Since by assumption of $\alpha \in G'(p) \subset S$, $(a_k\alpha)$ equidistributes in $[0,1]$, we thus obtain for sufficiently large $n$
\[({2^{n+1}})^{-\sigma}\#\{2^n \leq k \leq 2^{n+1}: a_k\alpha \in U'\}
\gg ({2^{n+1}})^{-\sigma}2^n \lambda(U) \gg \lambda(U)\sum_{2^n \leq k \leq 2 ^{n+1}} \lambda(B_k(\alpha)),
\]
proving \eqref{equid}. We now apply the Chung--Erd\H{o}s inequality to obtain
\[\lambda\left(\bigcup_{n_j \leq k,l < n_{j+1}} B_k \cap U \right)
\geq \frac{
\left(\sum_{n_j \leq k < n_{j+1}} \lambda(B_k \cap U)\right)^2
}{\sum_{n_j \leq k,l < n_{j+1}} \lambda(B_k \cap B_l \cap U)}
\geq \frac{
\lambda(U)^2\left(\sum_{n_j \leq k,l \leq j < n_{j+1}} \lambda(B_k)\right)^2
}{C\sum_{n_j \leq k,l < n_{j+1}} \lambda(B_k \cap B_l)} \geq \frac{\lambda(U)^2}{C'p}.
\]
Since this holds true for infinitely many $j$, we immediately deduce that 
\[\lambda(\limsup_{k \to \infty} B_k(\alpha) \cap U) \geq \frac{\lambda(U)^2}{C'p}.\]
An application of Proposition \ref{BDVdensitylemma} with $A = \limsup_{k \to \infty} B_k(\alpha)$ and $f(x) = \frac{x^2}{C'p}$ proves that 
\[\lambda(\limsup_{k \to \infty} B_k(\alpha)) = 1.\]

We now have shown that for $\alpha \in G'(p)$, $\lambda(\limsup_{k \to \infty} B_k(\alpha)) = 1$. 
Since $\mu(G'(p)) \geq 1 - 1/p$, the statement of the theorem follows by $p \to \infty$.

\section{Proof of Theorem~\ref{lebaefiber}}

The following lemma is standard. Since the proof is short, we nevertheless provide it for completeness.

\begin{lemma}\th\label{pairwiseind}
  The events
  \begin{equation*}
    A_n = \set*{(\alpha, \gamma)\in[0,1]^2 : \norm{a_n\alpha - \gamma}
      < \psi(n)} \subseteq [0,1]^2
  \end{equation*}
  are pairwise independent (with respect to $\lambda\otimes \lambda$).
\end{lemma}

\begin{proof}
  Let $m,n\in\NN$ with $m\neq n$. Let
  $\bone_m:= \bone_{(-\psi(m), \psi(m))}$ and
  $\bone_n:= \bone_{(-\psi(n), \psi(n))}$. Then
  \begin{align*}
    \meas(A_m\cap A_n)
    &= \int_{[0,1]^2} \bone_{A_m}\bone_{A_n} \, d(\lambda\otimes\lambda)\\
    &= \iint \bone_m(a_m\alpha - \gamma)\bone_n(a_n\alpha - \gamma)\, d\alpha\, d\gamma \\
    &= \iint \parens*{\sum_{k\in\ZZ}c_k^{(m)} e(k(a_m\alpha - \gamma))}\parens*{\sum_{\ell \in\ZZ}c_\ell^{(n)} e(\ell(a_n\alpha - \gamma))}\, d\alpha\, d\gamma 
  \end{align*}
  where $e(x):=\exp(2\pi i x)$ and $c_k^{(m)}$ and $c_\ell^{(n)}$ are the
  coefficients of the Fourier expansions of $\bone_m$ and $\bone_n$,
  respectively. Continuing, we have
  \begin{align*}
    \meas(A_m\cap A_n)
    &= \iint \sum_{k,\ell\in\ZZ}c_k^{(m)}c_\ell^{(n)} e(k(a_m\alpha - \gamma)) e(\ell(a_n\alpha - \gamma))\, d\alpha\, d\gamma  \\
    &= \sum_{k,\ell\in\ZZ}c_k^{(m)}c_\ell^{(n)} \parens*{\int e((ka_m+\ell a_n)\alpha )\, d\alpha} \parens*{\int e(-(k+\ell)\gamma)) \, d\gamma}.
  \end{align*}
  The last integral only contributes in the terms where $k+\ell = 0$, so
  \begin{equation*}
    \meas(A_m\cap A_n)
    = \sum_{k\in\ZZ}c_k^{(m)}c_{-k}^{(n)} \parens*{\int e(k(a_m- a_n)\alpha )\, d\alpha}.
  \end{equation*}
  Since $a_m\neq a_n$, the integral above only contributes when
  $k=0$, leaving
  \begin{equation*}
    \meas(A_m\cap A_n) = c_0^{(m)}c_0^{(n)} = \meas(A_m)\meas(A_n), 
  \end{equation*}
  confirming that the events $A_m$ and $A_n$ are independent.
\end{proof}

\begin{proof}[Proof of \th~\ref{lebaefiber}]
  By~\th~\ref{pairwiseind} and the second Borel--Cantelli lemma, the set
  \begin{equation*}
    W(\psi, \ba) :=\limsup_{n\to\infty} A_n
  \end{equation*}
  has full Lebesgue measure in $[0,1]^2$. By Fubini's theorem, almost
  every vertical fiber of $W(\psi, \ba)$ must also have full
  one-dimensional Lebesgue measure in $[0,1]$. In other words, 
  $\meas(W(\psi, \ba, \alpha))=1$ for almost every $\alpha\in[0,1]$.
\end{proof}

\section{Proof of Theorem~\ref{jarnik}}

\begin{proof}[Proof of~\th~\ref{jarnik}]
  Let $\ba = (a_n)$ and $\psi$ be as stated in the theorem. For each $n$, let
  \begin{equation*}
      B_n := B_n^1 = \set*{\gamma\in[0,1] : \norm{a_n\alpha - \gamma}
      < \psi(n)},
  \end{equation*}
  so that $W(\psi,\ba,\alpha) = \limsup B_n$. 
  
  We first dispatch the case where $\psi(n)\to 0$ does not hold. In this case, there is a sub-sequence $(a_{n_j})$ along which $\psi(n_j)$, and hence $\meas(B_{n_j})$, is bounded below by some $c>0$. Since for $\meas$-almost every $\alpha\in[0,1]$ the sequence $(a_{n_j}\alpha)$ equidistributes, it follows that for such $\alpha$ we have 
  \begin{equation*}
      \limsup_{j\to\infty} B_{n_j} = [0,1]
  \end{equation*}
  and therefore, $W(\psi, \ba, \alpha) = [0,1]$. Clearly, then, $\calH^f(W(\psi, \ba, \alpha)) = \calH^f([0,1])$ for any dimension function $f$. 
  
  We may now assume that $\psi(n)\to 0$. 
  Let $f$ be a dimension function such that
  $x^{-1} f(x)$ is monotonic. Suppose
  \begin{equation*}
    \sum_{n=1}^\infty f(\psi(n)) = \infty. 
  \end{equation*}
  Then, by~\th~\ref{lebaefiber}, for $\meas$-almost every $\alpha$ we
  have
  \begin{equation*}
    \meas(W(f\circ\psi, \ba, \alpha))=1,  
  \end{equation*}
  or, in other words,
  \begin{equation*}
    \calH^1(\limsup B_n^f) = \calH^1([0,1])
  \end{equation*}
  where
  \begin{equation*}
    B_n^f = \set*{\gamma\in[0,1] : \norm{a_n\alpha - \gamma}
      < f(\psi(n))}. 
  \end{equation*}
  By the mass transference principle (\th\ref{MTPthm}), we now have
  \begin{equation*}
    \calH^f(\limsup B_n^1) = \calH^f(W(\psi, \ba, \alpha)) = \calH^f([0,1]).
  \end{equation*}
  On the other hand, if 
  \begin{equation*}
    \sum_{n=1}^\infty f(\psi(n))  < \infty, 
  \end{equation*}
  then 
  \begin{equation*}
    \calH^f(\limsup B_n^1) = \calH^f(W(\psi, \ba, \alpha)) = 0
  \end{equation*}
  by the Hausdorff--Cantelli lemma.
\end{proof}

\begin{proof}[Proof of~\th\ref{jarnikbesicovitch}]
  Let $\psi_\sigma(n) = n^{-\sigma}$ and consider the dimension function
  $f(x)=x^{1/\sigma}$. Then we have
  \begin{equation*}
    \sum_{n=1}^\infty f(\psi_\sigma(n)) = \infty, 
  \end{equation*}
  and the divergence part of~\th~\ref{jarnik} implies that for
  $\meas$-almost every $\alpha$ we have
  \begin{equation*}
      \calH^{1/\sigma} (W(\psi_\sigma, \ba, \alpha)) = \infty,
    \end{equation*}
      hence $\dim_H(W(\psi_\sigma, \ba, \alpha)) \geq 1/\sigma$. On the other hand, if
  $\tau < \sigma$ and $g(x) = x^{1/\tau}$, then
  \begin{equation*}
    \sum_{n=1}^\infty g(\psi_\sigma(n)) < \infty, 
  \end{equation*}
  and \th\ref{jarnik} gives $\calH^{1/\tau}(W(\psi_\sigma,\ba,\alpha)) = 0$,
  hence $\dim_H(W(\psi_\sigma,\ba,\alpha)) \leq 1/\tau$. Since $\tau$ was
  arbitrary, we may conclude that given $\sigma\geq 1$, we have
  $\dim_H(W(\psi_\sigma, \ba, \alpha))=1/\sigma$ for almost every
  $\alpha\in [0,1]$.

  By intersecting countably many full-measure sets, we have that for
  $\meas$-almost every $\alpha\in[0,1]$, the formula
  \begin{equation*}
    \dim_H(W(\psi_\sigma, \ba, \alpha))=\frac{1}{\sigma}
  \end{equation*}
  holds for all $\sigma\in\QQ\cap [1, \infty)$. But in this case, it
  holds for all $\sigma\geq 1$, by the density of $\QQ$ and the
  monotonicity of Hausdorff dimensions.
\end{proof}

\bibliographystyle{plain}
\bibliography{bibliography.bib}




\end{document}